\RequirePackage{ifpdf}
\ifpdf 
\documentclass[pdftex]{sigma}
\else
\documentclass{sigma}
\fi

\newcommand{\pa}{\partial}
\newcommand{\D}{\Delta}
\newcommand{\R}{\mathbb R}
\newcommand{\N}{\mathbb N}

\begin{document}

\allowdisplaybreaks
	
\renewcommand{\PaperNumber}{121}

\FirstPageHeading

\renewcommand{\thefootnote}{$\star$}

\ArticleName{Conformal Powers of the Laplacian\\ via Stereographic
Projection\footnote{This paper is a
contribution to the Proceedings of the 2007 Midwest
Geometry Conference in honor of Thomas~P.\ Branson. The full collection is available at
\href{http://www.emis.de/journals/SIGMA/MGC2007.html}{http://www.emis.de/journals/SIGMA/MGC2007.html}}}

\ShortArticleName{Conformal Powers of the Laplacian via Stereographic
Projection}

\Author{C. Robin GRAHAM}

\Address{Department of Mathematics, University of Washington,\\
Box 354350, Seattle, WA 98195-4350, USA}

\Email{\href{mailto:robin@math.washington.edu}{robin@math.washington.edu}}
\URLaddress{\url{http://www.math.washington.edu/~robin/}}

\AuthorNameForHeading{C.R.~Graham}

\ArticleDates{Received November 17, 2007; Published online December 15, 2007}

\Abstract{A new derivation is given of Branson's factorization formula for
the conformally invariant operator on the sphere whose principal part is
the $k$-th power of the scalar Laplacian.  The derivation deduces Branson's
formula from knowledge of the correspon\-ding conformally invariant operator
on Euclidean space (the $k$-th power of the Euclidean Laplacian) via
conjugation by the stereographic projection mapping.}


\Classification{53B20}

\rightline{\it Dedicated to the memory of Tom Branson}

\section{Introduction}

The powers of the Laplacian on $\R^n$ satisfy an invariance property with
respect to conformal motions.  If $C$ is a conformal transformation
satisfying $C^*g_E = \Omega^2 g_E$, where $g_E$ denotes
the Euclidean metric and $\Omega$ is the conformal factor, then
\begin{gather}\label{translaw}
\Delta^{k} =
\left(C^{-1}\right)^*\Omega^{-n/2-k}\Delta^{k}\,\Omega^{n/2-k}C^*,
\qquad k\in \N,
\end{gather}
where the powers of $\Omega$ act by multiplication.
This observation is the motivation for
consideration of the ``conformally invariant powers of the Laplacian''
on a general curved conformal manifold (see~\cite{GJMS}).  In~\cite{B},
Tom Branson
derived the explicit form of such operators on the sphere $S^n$.  He showed
that any operator on $S^n$ which satisf\/ies the transformation
law analogous to~\eqref{translaw}, where now $C$ is a conformal
transformation of $S^n$ with conformal factor $\Omega$, necessarily is a
multiple of
\begin{gather}\label{sphere}
\prod_{j=1}^k(\D_{S}-c_j), \qquad c_j =
(\tfrac{n}{2}+j-1)(\tfrac{n}{2}-j).
\end{gather}
Here $\D_S$ denotes the Laplacian on the sphere, and our
sign convention is $\D=\sum \pa_i^2$ on $\R^n$.   To~prove this, he
introduced what are now called spectrum generating functions, by showing
how to use inf\/initesimal conformal invariance to derive the full spectral
decomposition of such an invariant operator from knowledge of its
eigenvalue on a single spherical harmonic.  Branson also used this argument
to give the
form of the pseudodif\/ferential intertwining operators satisfying a~transformation law analogous to \eqref{translaw} but involving more
general, possibly complex, powers of~$\Omega$.

There are now (at least) two other derivations of the
factorization \eqref{sphere}.  In \cite{Go}, \eqref{sphere} is derived via
a construction using tractors, and in \cite{FG} by explicit
solution of the algorithm of \cite{GJMS} in terms of the dual Hahn
polynomials, a family of discrete orthogonal polynomials.  Both of these
derivations
show that the same formula gives a conformally invariant operator
for any Einstein metric whose scalar curvature agrees with that of $S^n$.
This can also be deduced directly from Branson's result for $S^n$ and the
form of the GJMS algorithm; see the discussion in \cite{FG}.
A rescaling gives the corresponding formula for general Einstein metrics.

In this note we give a direct argument relating the operator $\Delta^k$
on $\R^n$ and the opera\-tor~\eqref{sphere} on $S^n$ under stereographic
projection.  Thus the conformal invariance of the operator
\eqref{sphere} is a~consequence of \eqref{translaw}.
The case $k=1$ is the Yamabe operator, whose conformal invariance, and
therefore whose behavior under stereographic projection, is
well-known.  The argument here deduces the relation
for $k>1$ from the case $k=1$ together with a calculation of pullback under
stereographic projection.
From this perspective, the constants $c_j$ for $j>1$ are manufactured from
$c_1$ by the stereographic projection mapping.

The derivation presented here is the analogue in the conformal case of an
argument in \cite{Gr} relating CR invariant operators on
odd-dimensional spheres to
corresponding operators on the Heisenberg group via the
Cayley transform.  The
CR case is more complicated: there is a 1-parameter
family of invariant operators for each $k$, and the operators on the
Heisenberg group are not powers of a f\/ixed operator, but rather are
products of various of the Folland--Stein operators.

\section{Derivation}

Let $\Phi:S^n\setminus \{p\} \rightarrow \R^n$ be stereographic
projection:
\[
\Phi(x',x_{n+1}) = x'(1+x_{n+1})^{-1}=y
\]
for $x'\in \R^n$ and $|x'|^2 + x_{n+1}^2 = 1$, where $p = (0,-1)$ is the
south pole.  One has
\[
\Phi^*\left(\frac{2}{1+|y|^2}\right)=1+x_{n+1}.
\]
The map $\Phi$ is conformal:
\[
\Phi^*g_E = (1+x_{n+1})^{-2} g_S.
\]
Def\/ine
$M^w:C^{\infty}(S^n\setminus \{p\})\rightarrow C^{\infty}(S^n\setminus
\{p\})$ by
\[
M^wf = (1+x_{n+1})^w f
\]
and $M_w:C^{\infty}(\R^n)\rightarrow C^{\infty}(\R^n)$ by
$M^w\Phi^* = \Phi^*M_w$, so that
\[
M_wf = 2^w(1+|y|^2)^{-w}f.
\]
The Yamabe operator on the sphere is
$Y=\D_S - c_1$, and its conformal invariance implies
\begin{gather}\label{yamabe}
YM^{1-n/2}\Phi^* = M^{-1-n/2}\Phi^*\D
\end{gather}
acting on functions on $\R^n$.

\begin{proposition}
For $k \in \mathbb N$,
\begin{gather}\label{main}
\left(\prod_{j=1}^k(\D_{S}-c_j)\right)M^{k-n/2}\Phi^* =
M^{-k-n/2}\Phi^*\D^k.
\end{gather}
\end{proposition}

The analogue of \eqref{translaw} for the operator \eqref{sphere} under
conformal
transformations of $S^n$ follows from \eqref{translaw} and \eqref{main},
since conjugation by $\Phi$ maps conformal transformations
of $\R^n$ to conformal transformations of $S^n$.

The proof begins by noting that $c_1-c_j=j(j-1)$, so that the left hand
side of \eqref{main} may be written as
\[
\left(\prod_{j=1}^k\left(Y+j(j-1)\right)\right)M^{k-n/2}\Phi^*.
\]
Now pass $\Phi^*$ through each term using \eqref{yamabe} and then cancel
the $\Phi^*$ to obtain that \eqref{main} is equivalent to the following
identity on $\R^n$:
\begin{gather}
[\D+k(k-1)M_2]M_{-2}[\D+ (k-1)(k-2)M_2]M_{-2}\cdots
[\D+2M_2]M_{-2}\D\nonumber\\
\qquad {}= M_{1-k}\D^k M_{1-k}.\label{Rn}
\end{gather}

The identity \eqref{Rn} can be proved by induction on $k$.  The induction uses
some commutator identities.  Denote by $X=\sum y_i\pa_{y_i}$ the Euler
vector f\/ield on $\R^n$.  The commutator identities are:
\begin{gather}\label{comm1}
[\D,X]=2\D,
\\
\label{comm2}
[X,M_w]=-w|y|^2M_{w+1},
\\
\label{comm3}
[\D,M_w] = -wM_{w}\left(2X +n-(w-1)M_1|y|^2\right)M_1,
\\
\label{comm4}
[\D^k,M_{-1}]=k\left(2X+n + 2(k-1)\right)\D^{k-1}.
\end{gather}
The f\/irst three are just direct calculations.  The last is an easy
induction on $k$.  Equation~\eqref{comm3} has been written in the form
above because
this is advantageous below, but it is easily seen using~\eqref{comm2} that
this may also be written perhaps a little more naturally as
\[
[\D,M_w] = -wM_{w+1}\left(2X +n-(w+1)M_1|y|^2\right).
\]
 In this form it is
clear that
the $k=1$ case of \eqref{comm4} is the $w=-1$ case of \eqref{comm3}.

Now prove \eqref{Rn} by induction.  The $k=1$ case is a tautology.
Assuming the result for $k$ and substituting this in the left hand side
for $k+1$ gives
\[
[\D+k(k+1)M_2]M_{-2}M_{1-k}\D^k M_{1-k},
\]
which equals
\begin{gather*}
\D M_{-k-1}\D^k M_{1-k} + k(k+1)M_{1-k}\D^k M_{1-k}\\
\qquad{} =M_{-k}\D M_{-1}\D^k M_{1-k} +
  [\D,M_{-k}]M_{-1}\D^kM_{1-k}+k(k+1)M_{1-k}\D^k M_{1-k}\\
\qquad {}=M_{-k}\D^{k+1}M_{-k} -M_{-k}\D[\D^k,M_{-1}]M_{1-k}
+[\D,M_{-k}]M_{-1}\D^kM_{1-k}\\
\qquad \phantom{={}}{} +k(k+1)M_{1-k}\D^k M_{1-k}.
\end{gather*}
Upon substituting \eqref{comm4} and \eqref{comm3} and then using
\eqref{comm1} to commute the $\D$ through the $X$ which arises in the
second term and f\/inally simplifying, one f\/inds that the last three terms
add up to 0, thus completing the induction step.

\subsection*{Acknowledgments}
This research was partially supported by NSF grant \# DMS
 0505701.

\pdfbookmark[1]{References}{ref}

\LastPageEnding

\end{document}